\documentclass[11pt,a4paper]{amsart}
\usepackage{amsthm}
\usepackage{amsfonts, amssymb, amsmath, amscd, latexsym}
\usepackage{psfrag, eucal, enumerate, latexsym}
\usepackage[twoside , dvips]{geometry}
\usepackage[all]{xy}

\usepackage[latin1]{inputenc}

\newtheorem{thm}{Theorem}[section]
\newtheorem{lem}[thm]{Lemma}

\theoremstyle{definition}
\newtheorem{rmk}[thm]{Remark}

\newtheorem{obs}[thm]{Observation}

\newcommand{\Z}{\mathbb Z} \newcommand{\C}{\mathbb C}
 
\newcommand{\p}{\mathbb P} \newcommand{\G}{\mathbb G}

\newcommand{\OO}{\mathcal{O}} \newcommand{\FF}{\mathcal{F}} \newcommand{\GG}{\mathcal{G}} \newcommand{\EE}{\mathcal{E}} \newcommand{\HHH}{\mathcal{H}} 

\DeclareMathOperator{\HH}{H}
\DeclareMathOperator{\hh}{h}
\DeclareMathOperator{\sat}{sat}
\DeclareMathOperator{\Ext}{Ext} \DeclareMathOperator{\Sym}{Sym}
\DeclareMathOperator{\Hom}{Hom} \DeclareMathOperator{\im}{Im}
\DeclareMathOperator{\rk}{rk}
\DeclareMathOperator{\norm}{norm}
\DeclareMathOperator{\even}{even} \DeclareMathOperator{\odd}{odd} \DeclareMathOperator{\ch}{char}
\DeclareMathOperator*{\Oplus}{\oplus}

\SelectTips{cm}{}
\newdir{ >}{{}*!/-5pt/@{>}}

\begin{document}
\sloppy


\title{A geometric construction of Tango bundle on $\p^5$}

\author{Daniele Faenzi}
\email{faenzi@math.unifi.it}
\address{Dipartimento di Matematica ``U.~Dini'', Universit\`a di
 Firenze, Viale Morgagni 67/a, I-50134, Florence, Italy}
\urladdr{http://www.math.unifi.it/\~{}faenzi/}
\thanks{{\em 2000 Mathematics Subject Classification} Primary 14F05; Secondary 14-04,14J60}
\thanks{The author was partially supported by Italian Ministery funds
and the EAGER contract HPRN-CT-2000-00099.}

\begin{abstract} 
The Tango bundle $T$ over $\p ^5$ is proved to be the pull-back of the
twisted Cayley bundle $C(1)$ via a
map $f \colon \p ^5 \rightarrow Q_5$ existing only in characteristic
2. The Frobenius morphism $\varphi$ factorizes via such $f$. \\
Using $f$ the cohomology of $T$ is computed in terms of 
$S \otimes C$, $\varphi^*(C)$, $\Sym^2(C)$ and $C$, while these are
computed by applying Borel-Bott-Weil theorem. \\ 
By machine-aided computation the mimimal resolutions of $C$ and $T$
are given; incidentally the matrix presenting the spinor bundle $S$ over
$Q_5$ is shown.

\end{abstract}


\maketitle

\section{Introduction}

The well-known Hartshorne conjecture states, in particular, that there
are no indecomposable rank-2 vector bundles on $\p^n$, when $n$ is
greater than 5. However, one of the few rank-2 bundles on $\p^5$ up to
twist and pull-back by finite morphisms is the Tango bundle $T$ first
given in \cite{tan-morphisms-I}. Later Horrocks in \cite{horrocks} and Decker Manolache
and Schreyer in \cite{decker} descovered that it can be obtained starting from Horrocks rank-3
bundle: anyway it only exists in characteristic 2.

Here we prove that $T$ is the pull-back of the
twisted Cayley bundle $C(1)$ (defined over any field) via a
map $f \colon \p ^5 \rightarrow Q_5$ existing only in characteristic
2. This allows to
compute its cohomology by Borel-Bott-Weil theorem (see \ref{bott-cohomology}). We observe non-standard cohomology of $\Sym ^2C$
in characteristic 2, although $C \otimes S$ behaves standardly
(i.e. as if $\ch (k)=0$). 

Different ways to prove that Tango's equations
actually give $C(1)$ are shown in section (\ref{sec-2}). Lastly in
section (\ref{sec-3}) we give  the minimal resolution of $T$ and some
more computational remarks. 
We make extensive use of {\tt Macaulay2} computer
algebra system: computation-related material can be found on the url:
{\tt http://www.math.unifi.it/\~{}faenzi/tango/} \\
I would like to thank Wolfram Decker for the observations 
on the Horrocks bundle and his help on machine computation and Edoardo
Ballico for his careful
remarks about Ekedahl's work and representations in positive characteristic. Also I am indebted with Giorgio
Ottaviani, who posed to me this problem.
\vspace{0.3cm}

So let
        $$Q_5 = \{ z_0^2+ z_1z_2+z_3z_4+z_5z_6 = 0 \} \subset \p^6$$
be the 5-dimensional smooth quadric over an algebraically closed field $k$. We denote by $\xi$ (respectively $\eta$, $\zeta$) the generator of
$A^1(\p ^5)$ (respectively of $A^1(Q_5)$, $A^3(Q_5)$),
so that
        \begin{align*}
         A(\p^5)   =      \Z[\xi] / (\xi ^6) \quad  A(Q_5)       =      \Z[\eta , \zeta] / (\eta ^3 - 2 \zeta, \eta^6)
        \end{align*}
On the
coordinate ring $S(\p^5)$  we use variables $x_i$'s while on
$S(Q_5)$ we use $z_j$'s.

\section{The bundle on $\p ^5$} \label{sec-1}

Let $k$ be an algebraically closed field.
On $Q_5=G_2 / P(\alpha_1) $ we have the Cayley bundle $C$, coming from the
standard representation of the semisimple part of the parabolic group
$P(\alpha_1)$, where $\alpha_1$ is the shortest root in the Lie
algebra of $G_2$, the exceptional Lie group.  $C$ 
is irreducible $G_2-$homogeneous with maximal weight $\lambda_2 -
2\lambda_1$. $C(2)$ (weight $\lambda _2$) is globally generated
and $\hh^0(C(2))=14$.\\
$C$ is the cohomology of a monad:
        $$\OO(-1) \longrightarrow S  \longrightarrow \OO $$
where $S$ is the spinor bundle. $C$ has rank 2 and Chern classes $(-1,1)$.  The only
non-vanishing intermediate cohomology groups are $\HH^1(C)=\HH^4(C(-4))
=k$. All this is done in \cite{ott-cayley} and follows easily from
\cite[Proposition 5.4]{jantzen} in any characteristic. \\
When $\ch (k)=2$ we have a map $f \colon \p ^5
\rightarrow Q_5$ having the expression
        $$f : (x_0:\ldots:x_5) \mapsto (x_0x_1 + x_2x_3 + x_4x_5,x_0^2,x_1^2,x_2^2,x_3^2,x_4^2,x_5^2)$$
and, letting $\varphi$ be the Frobenius and $\pi$ the projection from
$(1:0:0:0:0:0:0)$ we have:
        \begin{equation}\label{pull-back}
        \xymatrix @C+3ex{
                                                               & Q_5 \ar^-{\pi}[d]      & \varphi^*(\OO_{\p^5}(1))     =  \OO_{\p^5}(2) \\
        \p^5 \ar^-{f}[ur] \ar_-{\varphi}[r]                    & \p^5                   & \pi^*(\xi) =  \eta \qquad  \pi^*(\xi^2) = \eta ^2 \qquad \pi^*(\xi^3) = \eta ^3 =  2 \zeta 
        }
        \end{equation}  
We define  $T = f^*(C(1))$. The rank$-2$ vector bundle $T$ is the main subject of this paper. 

\vspace{0.3cm}
By Hirzebruch-Riemann-Roch
theorem we have the following easy observation
\begin{obs}
$T$ has Chern classes $(2,4)$ and
        $$\chi(T(t)) = \frac{1}{60} {t}^{5}+\frac{1}{3} {t}^{4}+\frac{25}{
         12}  {t}^{3}+\frac{11}{ 3} {t}^{2}-\frac{51}{ 10} t-14$$
\end{obs}
The proof of the following theorem on the cohomology of $T$ is left at
the end of this section.
\begin{thm} \label{bott-cohomology}
The non-zero cohomology of $T$ is the following:
        \begin{align*}
        & \hh^0(T(t)) =  \chi(T(t)) \qquad \mbox{for $t\geq 2$} \\
        & \hh^1(T(-2)) =  1 \quad \hh^1(T(-1)) = 7 \quad \hh^1(T) = 14 \quad  \hh^1(T(1)) = 13 \quad \hh^1(T(2)) = 1 \quad \\
        & \hh^2(T(-3)) = 1
        \end{align*}
and their Serre-dual $\hh^{5-i}(T(t)) = \hh^i(T(-t-8))$. 
$T'=T(-4)$ is the only twist
for which all cohomology groups vanish.
\end{thm}

For this we first need an observation on the map $f$. 
\begin{obs}\label{fstar}
For the map $f$ above we have 
        \begin{align}
        & f_{*}(\OO_{\p^5})  =  \OO \oplus \OO(-1)^{14} \oplus \OO(-2)
        \nonumber \\
        \label{fstar1}
        & f_{*}(\OO_{\p^5}(1))  = S \oplus \OO^6 \oplus \OO(-1)^6
        \end{align}
\end{obs} 
\begin{proof}
Let $\FF = f_{*}(\OO_{\p^5}), \GG = f_{*}(\OO_{\p^5}(1))$. The map $f$ is a $16\colon 1$ cover,
because the Frobenius is $32 \colon 1$ and the projection $\pi$ is
$2\colon1$. Then $\FF$ and $\GG$ are rank-16 vector bundles, whose 
cohomology one can read from the Leray degenerate spectral sequence. Indeed since $R^i(f_*) = 0$ (for $i>0$) we have
        \begin{align*}
        & \HH^i(Q_5,\FF(t)) = \HH^i ( \p^5, f^* (\OO_{Q_5}(t)) = \HH^i (\p^5 , \OO_{\p^5} (2t)) \\
        & \HH^i(Q_5,\GG(t)) = \HH^i ( \p^5, f^* (\OO_{Q_5}(1)(t)) =  \HH^i (\p^5 , \OO_{\p^5} (2t+1)) 
        \end{align*}
for $0\leq i\leq 5$ and every $t$. This says that $\FF$ and $\GG$ have no
intermediate cohomology, hence by \cite{kapranov} or \cite{buchweitz-greuel-schreyer}
they must decompose as sum of spinors $S$ and line bundles, up to
twist (although actually Kapranov's setting is over $\C$). For $\GG$ this implies, by a computation on
the Euler characteristic, that the only choice is the one stated while
for $\FF$ we have a priori two possibilities: the one stated above or
$S \oplus S(-1) \oplus \OO \oplus  \OO(-1)^6 \oplus \OO(-2)$ \\
Now the above formula says that the polynomial ring $S(\p^5)$
decomposes as module over $S(Q_5)$ (under the action given by $f$) as $S(\p^5)_{\even} \oplus
S(\p^5)_{\odd}$ where:
        $$S(\p^5)_{\even} = \Oplus_{t \in \Z} \HH^0(\p^5,\OO(2t)) \qquad 
        S(\p^5)_{\odd} = \Oplus_{t \in \Z} \HH^0(\p^5,\OO(2t+1))$$
For $\FF$ we have to compute explicitely a presentation of the $S(Q_5)-$module
$S(\p^5)_{\even}$. We need $e_0$ to generate $k=S(\p^5)_{0}$ and
$e_{ij}$ to generate the monomial $x_ix_j$ ($i\neq j$) in
$S(\p^5)_{2}$, thus obtaing a map
        $$\Phi \colon S(Q_5) \oplus S(Q_5)(-1) ^{15}
        \longrightarrow S(\p^5)_0 \oplus S(\p^5)_{2}$$
But the coordinate $e_{45}$ is redundant, since $z_0 \Phi(e_0) +
\Phi(e_{01}) + \Phi(e_{23}) = \Phi(e_{45})$. Now for $S(\p^5)_{4}$:
the terms containing $x_i^2$ already lie in the image (got by the
action of $z_{i-1}$), and in fact we just have to fix $x_0x_1x_2x_3$
because, e.g. $x_0x_1x_2x_4 = z_3\Phi(e_{34})+z_5\Phi(e_{25})+z_0\Phi(e_{24})$ and 
$x_0x_1x_4x_5=x_0x_1x_2x_3+z_1z_2\Phi(e_0)+z_0\Phi(e_{01})$. Thus we get a generator in degree $2$ (and
no syzygy); moreover $x_0x_1x_2x_3x_4x_5 =
z_1z_2\Phi(e_{23})+z_3z_4\Phi(e_{01})+z_0(x_0x_1x_2x_3)$, so that $S(\p^5)_{6}$, and in
fact all $S(\p^5)_{\even}$, is also covered.

The presentation could be computed also for $S(\p^5)_{\odd}$, where
one finds as syzygy of a map $S(Q_5)^{6} \oplus S(Q_5)(-1)^{14} \rightarrow S(\p^5)_{\odd}$ the matrix
giving the spinor bundle described later in (\ref{spinmatrix}), thus getting again (\ref{fstar1}).
\end{proof}

\begin{rmk}

It seems likely that the observation (\ref{fstar}) can be extended to any odd
dimension. Here we only mention that for $f \colon \p^3
\longrightarrow Q_3$ we get $f_*(\OO_{\p^3}) = S \oplus \OO_{Q_3}
\oplus \OO_{Q_3} (-1)$ and $f_*(\OO_{\p^3}(1))= \OO_{Q_3}^4$
which can be computed by the presentation or by Euler characteristic.

Lastly, one may notice that $\pi_*(S)=\OO_{\p^5}(-1)^8$ and clearly
 $\pi_*(\OO_{\p^5}(t))=\OO(t) \oplus \OO(t-1)$ so that the
extension $0 \rightarrow S(-1) \rightarrow \OO(-1)^8 \rightarrow S
\rightarrow 0$ becomes split after $\pi_*$ (actually this holds in any characteristic). This agrees with the
formula
        \begin{align*}
        & \varphi_*(\OO_{\p^5}) = \OO \oplus \OO(-1)^{15} \oplus \OO(-2)^{15} \oplus \OO(-3) \\
        & \varphi_*(\OO_{\p^5}(1)) = \OO^6 \oplus \OO(-1)^{20} \oplus \OO(-2)^{6}
        \end{align*}
\end{rmk}

Next we compute the cohomology of $\Sym ^2 C$, $C \otimes C$, $S
\otimes C$. First notice that if $V$ is the $SL(2)-$representation
giving $C$, when $\ch(k)=2$ the representation $\Sym ^2 V$ 
(having weight $2\lambda _2 - 4\lambda _1 $) will not be irreducible
 (recall that in finite characteristic $SL(2)$ is not linearly
reductive, check \cite{nagata}), on the contrary, letting $C^{[2]} = \varphi ^* (C)$, we have the
non-split exact sequence
        \begin{equation} \label{C2}
        0 \longrightarrow C^{[2]} \longrightarrow \Sym ^2 C 
                \longrightarrow \OO(-1) \longrightarrow 0
        \end{equation}
We also have
        \begin{align} \label{phiQ5}
        & \varphi_*(\OO_{Q_5}) = \OO \oplus \OO(-1)^{20} \oplus \OO(-2)^{7} \oplus S(-1) \\
        & \varphi_*(\OO_{Q_5}(1)) = \OO^7 \oplus \OO(-1)^{20} \oplus \OO(-2) \oplus S \nonumber
        \end{align}
Now again by Borel-Bott-Weil theorem (\cite[Proposition 5.4]{jantzen}) we 
know $\hh^i(\Sym^2 C(3)) = 0$, $\Sym(4)$ is globally generated and:
        \begin{align*}
        & \hh^0(\Sym(t)) = \chi(\Sym ^2 C(t)) =
\frac{1}{20}\,{t}^{5}+\frac{3}{8}\,{t}^{4}-{\frac
{27}{8}}\,{t}^{2}-{\frac{81}{20}}\,t \quad \text{for $t\geq 4$} \\
        &  \quad \hh^1(\Sym^2 C(2)) = 14        \quad \hh^1(\Sym^2 C(1)) = 7
        \quad \hh^0(\Sym^2 C(-1)) = 1
        \end{align*}
where in the above twists these are the only non-vanishing $\HH^i$'s. 
By Serre duality we only miss $\HH^i(\Sym ^2
C)$. Now since $R^i( \varphi _*) = 0$, by Leray spectral sequence,
(\ref{C2}) and (\ref{phiQ5}) we get:
        \begin{equation*}\begin{split}
              \HH^2(\Sym^2 C) = & \HH^2(C ^{[2]})
\overset{\text{Leray}}{=} \HH^2(C \otimes S(-1)) \overset{\text{Leray}}{=} \\
                & = \HH^2(C ^{[2]}(-1)) = \HH^2(\Sym^2 C(-1)) \overset{\text{Bott}}{=} k
        \end{split}\end{equation*}
and by the same argument $\hh^1(\Sym^2 C) = 1$, while the remaining
$\hh^i$ are zero. The same procedure yields the values of $\hh^i(S
\otimes C (t))$
        \begin{align} \label{cohomology-SC}
         & \hh^1(S\otimes C )=1 & \hh^1(S\otimes C (1))=6 \\
         & \hh^2(S\otimes C (-1))=1 & \hh^3(S\otimes C (-2))=1 \nonumber\\
         & \hh^4(S\otimes C (-4))=6 & \hh^4(S\otimes C (-5))=1 \nonumber
        \end{align}
Now we can prove theorem (\ref{bott-cohomology})
\begin{proof}[Proof of theorem (\ref{bott-cohomology})]
Since $R^i(f_*) = 0$ for $i>0$ we use the Leray degenerate spectral sequence that here reads:
        \begin{align} \label{leray}
        \HH ^i(\p^5,T(2t)) = & \HH^i(Q_5,f_*(\OO_{\p^5}) \otimes C(1+t)) \nonumber \\
        \HH ^i(\p^5,T(2t+1)) = & \HH^i(Q_5,f_*(\OO_{\p^5}(1)) \otimes C(1+t)) 
        \end{align}
The even part gives $\HH ^i(\p^5,T(2t))=\HH^i(Q_5, C(1+t)) \oplus
\HH^i(Q_5,C(t))^{14} \oplus \HH^i(Q_5,C(t-1))$, while the odd part gives
        $\HH^i(\p^5,T(2t+1)) = \HH^i(Q_5, C(1+t))^6 \oplus \HH^i(Q_5,
        C(t))^6 \oplus \HH^i(Q_5, S \otimes C(1+t))$ and these are known by the
        above formulas.
\end{proof}

Now we can write the Beilinson table of the normalized $T(-1)$ (rather
than of $T$) that reads:
        $$\left(\begin{array}{cccccc}
  0     & 0     & 0     & 0     & 0     & 0     \\
  1     & 0     & 0     & 0     & 0     & 0     \\
  0     & 1     & 0     & 0     & 0     & 0     \\
  0     & 0     & 0     & 1     & 0     & 0     \\
  0     & 0     & 0     & 0     & 1     & 7     \\
  0     & 0     & 0     & 0     & 0     & 0 
        \end{array}     \right)$$
Then we write $T(-1)$ as cohomology of a monad:
        $$ \OO(-1) \oplus \Omega ^4 (4) \longrightarrow \Omega ^2 (2)
        \oplus \Omega ^1 (1) \longrightarrow  \OO ^7 $$

\vspace{0.3 cm}
The following remark shows that our situation is more ``rigid'' than
one may expect.
\begin{rmk}
We notice that we have the diagram:
$$\xymatrix@C+2ex{
          & Q_5  \ar[d]_-{\pi} \ar^-{\varphi}[r]   & Q_5\\
          \p^5 \ar_{\varphi}[r] \ar^f[ur] & \p^5 \ar_f[ur] & 
        }$$
Then, as Edoardo Ballico pointed out to us, we may relate this
 framework to a description
given by Ekedahl in \cite[proposition 2.5]{eke-foliations} of finite morfisms $\psi$ (with
$Y$ smooth) that
factor through the Frobenius $\varphi$ in characteristic $p>0$:
        $$\xymatrix{
          Y  \ar[d]_-{\psi} \ar^-{\varphi}[r]   & Y\\
          \p^n \ar[ur] & 
        }$$
Ekedahl shows that in such a situation $Y$ is $Q_n$, $n$ is odd, $\psi$
is the projection from a point external to $Q_n$ and the
characteristic is 2. That is, precisely our setup.
\end{rmk}
\begin{rmk} The values $\hh^2(\Sym^2 C)=\hh^1(\Sym^2 C)=1$ exhibit
non-standard cohomology for the representation $\Sym^2 V$. Indeed
$3\lambda_2 - 3 \lambda_1$ is singular ($(3\lambda_2 - 3 \lambda_1,3
\alpha_1 + \alpha _2)=0$) so standard Borel-Bott-Weil theorem (i.e. in
characteristic 0) would give $\hh^i(\Sym^2 C)=0$. \\
Of course, we would have no such sequence as (\ref{C2}). Still, by
tensoring the monad defining $C$ by $C(t)$ we get
        \begin{gather*}
          0 \longrightarrow G^{\vee} \otimes C(t) \longrightarrow S \otimes
        C(t) \longrightarrow C(t) \longrightarrow 0 
        \\
          0 \longrightarrow C(t-1) \longrightarrow G^{\vee} \otimes C(t)
        \longrightarrow C \otimes C(t)  \longrightarrow 0
        \end{gather*}
whence we derive the values of $\hh^i(S \otimes C(t))$ from those of
$\Sym^2 C$, since if $\ch (k) \neq 2$ $C \otimes C = \Sym^2 C \oplus
\OO (-1) $. This way we would get the
same values as in (\ref{cohomology-SC}) if $\ch (k)=0$. \\
Finally, one can compute on {\tt Macaulay2} the values of the cohomology
of $C^{[2]}$ and check the correcteness of the above result.
\end{rmk}

\section{The Cayley bundle and Tango's equations} \label{sec-2}

In this sections we work over any algebraically closed field $k$ and prove that Tango's
equations in \cite{tan-morphisms-I} give $C(1)$. \\
First we complete Tango's $3 \times 6$ matrix by stacking 9 rows to get $A$ having everywhere rank 3 over the quotient ring $S(Q_4) = 
k[z_0,\ldots,z_6]/(z_0^2+z_1z_2+z_3z_4+z_5z_6)$
        \begin{equation} \label{matrix-A}
        A = {\scriptsize
        \left(
        \begin{array}{cccccc}
         z_0^2                  & 0                      & 0                      & z_1^2                 & z_1z_3+z_0z_6          & -z_0z_4+z_1z_5      \\
         0                      & z_0^2                  & 0                      & z_1z_3-z_0z_6         & z_3^2                  & z_0z_2+z_3z_5       \\
         0                      & 0                      & z_0^2                  & z_0z_4+z_1z_5         & -z_0z_2+z_3z_5         & z_5^2               \\
         0                      & z_0z_2-z_3z_5          & z_3^2                  & -z_0z_3-z_2z_6        & 0                      & z_2^2               \\
         z_5^2                  & 0                      & z_0z_4-z_1z_5          & z_4^2                 & -z_2z_4-z_0z_5         & 0                   \\
         -z_3^2                 & z_1z_3+z_0z_6          & 0                      & -z_6^2                & 0                      & -z_0z_3+z_2z_6      \\
         0                      & z_5^2                  & -z_0z_2-z_3z_5         & -z_2z_4+z_0z_5        & z_2^2                  & 0                   \\
         z_1z_3-z_0z_6          & -z_1^2                 & 0                      & 0                     & -z_6^2                 & z_0z_1+z_4z_6       \\
         z_0z_4+z_1z_5          & 0                      & -z_1^2                 & 0                     & -z_0z_1+z_4z_6         & -z_4^2              \\
         -z_2^2                 & -z_2z_4-z_0z_5         & z_0z_3-z_2z_6          & z_0^2                 & 0                           & 0              \\
         z_2z_4-z_0z_5          & z_4^2                  & z_0z_1+z_4z_6          & 0                     & z_0^2                  &0                    \\
         -z_0z_3-z_2z_6         & z_0z_1-z_4z_6          & -z_6^2                 & 0                     & 0                      & -z_0^2
        \end{array}
        \right)}
        \end{equation}
indeed we have:
        \begin{align*}
        & \sat (A) = \sat ((\wedge ^2 A))  = \sat ((\wedge ^3 A)) = (1) \\
        & \wedge ^4 A             = (0)
        \end{align*}

Now let $H$ be the sheaf $(\im (A) \widetilde{\; \;})^{\vee}(-1)$. Tango's bundle
$F$ is defined as the cokernel $0 \rightarrow \OO \rightarrow H 
\rightarrow F \rightarrow 0$.

One can ask {\tt Macaulay2} for the cohomology groups $\HH^i(H)$ and
the result is:
        \begin{align*}
        & \hh^i(Q_5,H(t))  =  0 \qquad \forall \hspace{1mm} t, \hspace{3mm} 0<i<5, \mbox{ except:} \\
        & \hh^1(Q_5,H(-1))  =  1
        \end{align*}
Therefore we have:
        $$\Ext^1(\OO,H(-1)) = k$$
the associated extension is:
        $$ 0 \longrightarrow H(-1) \longrightarrow W \longrightarrow \OO \longrightarrow 0 $$
where $W$ is a 4-bundle whose intermediate cohomology is forced to be
zero. Then by \cite{buchweitz-greuel-schreyer} or \cite{kapranov} and
Euler characteristic we get $W = S$ and $H(-1)=G^{\vee}$ so that $F =
C(1)$.

\vspace{0.3 cm} 
A different method would be proving that $C(1)$ and $F$ have the same
resolution. The resolution of $C(1)$ can be obtained on the computer by 
 in the following way. We first need a matrix for the spinor
bundle $S$ and we derive it from \cite{beauville-determinantal}. \\
We know that $S$ is a rank$-4$ bundle with no intermediate cohomology,
$S^{\vee}= S(1)$, and $\hh^0(S(1))=\hh^5(S(-5))=8$. \\
Then by Beilinson theorem $S(1)$ extended by zero to $\p^6$ has the resolution:
        $$ 0 \longrightarrow \OO_{\p^6}(-1)^8 \overset{B}{\longrightarrow} \OO_{\p^6}^8
        \longrightarrow S(1) \longrightarrow 0 $$
where now $B$, by the observations in
\cite{beauville-determinantal}, is an antisymmetric matrix whose determinant is the equation of
the quadric, to the power 4. This is done by the matrix
        \begin{equation} \label{spinmatrix}
        B = 
 \begin{pmatrix}0&
      0&
      0&
      {-{{z}}_{{3}}}&
      0&
      {-{{z}}_{1}}&
      {{z}}_{{5}}&
      {-{{z}}_{0}}\\
      0&
      0&
      {{z}}_{{3}}&
      0&
      {{z}}_{1}&
      0&
      {-{{z}}_{0}}&
      {-{{z}}_{{6}}}\\
      0&
      {-{{z}}_{{3}}}&
      0&
      0&
      {-{{z}}_{{5}}}&
      {{z}}_{0}&
      0&
      {-{{z}}_{{2}}}\\
      {{z}}_{{3}}&
      0&
      0&
      0&
      {{z}}_{0}&
      {{z}}_{{6}}&
      {{z}}_{{2}}&
      0\\
      0&
      {-{{z}}_{1}}&
      {{z}}_{{5}}&
      {-{{z}}_{0}}&
      0&
      0&
      0&
      {{z}}_{{4}}\\
      {{z}}_{1}&
      0&
      {-{{z}}_{0}}&
      {-{{z}}_{{6}}}&
      0&
      0&
      {-{{z}}_{{4}}}&
      0\\
      {-{{z}}_{{5}}}&
      {{z}}_{0}&
      0&
      {-{{z}}_{{2}}}&
      0&
      {{z}}_{{4}}&
      0&
      0\\
      {{z}}_{0}&
      {{z}}_{{6}}&
      {{z}}_{{2}}&
      0&
      {-{{z}}_{{4}}}&
      0&
      0&
      0\\
      \end{pmatrix}
        \end{equation}


Incidentally, we mention that by matrices written in such a fashion one
can obtain the spinor bundles over the quadric of any dimension.
Here by a standard mapping cone construction we get an infinite 2-periodic resolution of 
$C(1)$ of the form:
        $$ \textstyle {R^{\bullet}  
                        \overset{\delta}{\longrightarrow}     \OO_{Q_5}(-4)^{55} 
                        \rightarrow     \OO_{Q_5}(-3)^{49} 
                        \rightarrow     \OO_{Q_5}(-2)^{34} 
                        \rightarrow     \OO_{Q_5}(-1)^{14} 
                        \rightarrow     C(1) \rightarrow   0}
        $$
where 
        $$\textstyle{R^{\bullet} =
        \cdots  \overset{d}{\longrightarrow}    \OO_{Q_5}(-i-3)^{56}
                \overset{e}{\longrightarrow}    \OO_{Q_5}(-i-2)^{56} 
                \overset{d}{\longrightarrow}    \OO_{Q_5}(-i-1)^{56} 
                \overset{e}{\longrightarrow} 
         \cdots}$$
and this coincides with what we get by Tango's original
construction. 

\vspace{0.3cm}
We remark that the periodicity here is a standard behaviour, by \cite{eisenbud}. \\
The rank of the kernel of $\delta$ is 28 and again it must have no
intermediate cohomology. Then Euler characteristic shows that
it must be $S(-5)^{\oplus 7}$, i.e. the resolution actually reads:
        $$ \textstyle {0 \rightarrow    S(-5)^{\oplus 7}
                        \rightarrow     \OO_{Q_5}(-4)^{55} 
                        \rightarrow     \OO_{Q_5}(-3)^{49} 
                        \rightarrow     \OO_{Q_5}(-2)^{34} 
                        \rightarrow     \OO_{Q_5}(-1)^{14} 
                        \rightarrow     C(1) \rightarrow   0}
        $$
This agrees with \cite[pag. 197]{ottaviani-szurek} (here we consider
the dual spectral sequence) i.e. the above is
a {\em Kapranov sequence} for $C(1)$.

\vspace{0.3 cm}
Yet another method is the following. First one proves  the analogous for $Q_5$
of \cite[lemma 1]{tan-morphisms-I}, i.e.

\begin{lem}\label{enumerative} Let  $\rho \colon Q_n \longrightarrow \G(\p^k,\p^n)$ be a
 non-constant morphism, with $n$ odd and $k$ even, and let $\EE$ be the
 pull-back on $Q_n$ of the dual universal sub-bundle $U^{\vee}$ on $\G(\p^k,\p^n)$. 
 Then $k=\frac{n-1}{2}$, and:
        $$c_i (\EE)= 2a^i$$
for some positive integer $a$.
\end{lem}

Here $\rho^*(U)(1)=H$ defined above, where $\rho$ is given by
the matrix $A$ (\ref{matrix-A}): a row of $A(x)$ is a point in $\p ^6$
and $\rk (A(x))=3$, hence $A(x)$ represents a $\p^2 \subset \p^6$. \\
But we know that $\rho^*(U)$ contains the
sub-line-bundle $\OO(-1)$, so that
$a=1$.  So $H_{\norm} = H(-1)$. Then {\tt Macaulay2} gives
$\HH^0(H_{\norm}) = \HH^0(\wedge^2 H_{\norm}) =0$ thus $H$ is stable
by Hoppe's criterion. Hence we can conclude by \cite{ott-spinor} that
$\rho^*(U)=G^{\vee}$ and $F=C(1)$. 

\begin{proof}[Proof of lemma (\ref{enumerative})]
The proof is almost identical to the case of $\p^n$,
the only difference being 
that we have to work in $H^*(Q_n)$, where
$\eta^{\frac{n+1}{2}}=2\zeta$. \\
We can suppose $k$ even and $k\leq\frac{n-1}{2}$ because
$\G(\p^k,\p^n) \cong \G(\p^{n-k-1},\p^n)$ and we put
$a_i = c_i(\rho^*(\EE))$, $b_i = c_i(\rho ^*(Q))$ where $Q$ is the
universal quotient bundle. We have, in the ring $A(Q_5)[t]$ the
relation on Chern polynomials
        \begin{equation} \label{chern-1}
        c _{\EE^{\vee}} (t) \cdot c_{\rho^*(Q)} (t) = 1  
        \end{equation}
Now we can think of the coefficients in (\ref{chern-1}) as integers times some $\eta^r$, taking care
to replace $\zeta$ by $\frac{1}{2}\eta^{\frac{n+1}{2}}$, that is,
replacing $a_i$ (and $b_i$)
by $a_i'=\frac{1}{2}a_i$ (by $b_i'=\frac{1}{2}b_i$) whenever $i\geq\frac{n+1}{2}$. \\
Then  one proceeds exactly as in \cite[lemma 1]{tan-morphisms-I} and
\cite[lemma 3.3]{tan-projective}, and finds:
         \begin{align*}
         & k = \frac{n-1}{2} \\
         & a'_i = 2 a^i \qquad \mbox{for $i=1,\ldots,\frac{n-1}{2}=k$} \\
         & a'_{k+1} = a^{k+1}
         \end{align*}
and so $a_i = 2a^i$, for all $i$'s, as only for $i=k+1$ we have to
substitute $a_{k+1} = 2a'_{k+1}$.
\end{proof}

\section{further remarks} \label{sec-3}

Now we turn back to $\p^5$ and $\ch(k)=2$. As we have the equations
for $T$, {\tt Macaulay2} provides  the following resolution
        $$
        \xymatrix@C-2ex{
        0 \ar[r]
        & *\txt{$\OO^7(-7)$ \\ $\oplus$ \\ $\OO(-8)$} \ar[r] 
        & \OO^{49}(-6) \ar[r] 
        & \OO^{98}(-5) \ar[r]
        & \OO^{76}(-4)\ar[r]
        & *\txt{$\OO^{7}(-3)$ \\ $ \oplus $ \\ $\OO^{14}(-2)$} \ar[r]
        & T \ar[r] & 0
        }$$
One may check that this gives back the Chern classes $2,4$. \\
As Wolfram Decker pointed out to me, another way to get Tango's bundle
is by Horrocks bundle $\HHH$ in characteristic 2. Concretely, Horrocks 
becomes a non-split extension $0 \rightarrow T(-1) \rightarrow \HHH
        \rightarrow \OO \rightarrow 0$. \\
This allows to compute the cohomology of $\HHH$ in terms of $T$. 

Moreover in \cite{decker} one finds an explicit description of the maps in the Beilinson
monad. This provides also a check of computations. Denoting with $e_0,\ldots,e_5$ the canical basis of $E$ (the exterior algebra
over $V$), and using the natural isomorphism
$\Hom(\Omega^i(i),\Omega^j(j)) = \wedge ^{i-j} V = E_{j-i}$, $T$ is
the cohomology of the maps $\alpha$ and $\beta$:
        $$\begin{array}{c}
        \beta   = \left( \begin{array}{cc}
        e_0     &       e_4e_5               \\
        e_1     &       e_3e_5               \\
        e_2     &       e_3e_4               \\
        e_3     &       e_1e_2               \\
        e_4     &       e_0e_2               \\
        e_5     &       e_0e_1               \\
        0       &       e_0e_3+e_1e_4+e_2e_5 
        \end{array}\right)
        \\  \\
        \alpha  = \left( \begin{array}{ccc}
         e_0e_1e_2+e_3e_4e_5    & e_0e_1e_3e_4+e_0e_2e_3e_5+e_1e_2e_4e_5        & 0     \\
        e_0e_3+e_1e_4+e_2e_5    & e_0e_1e_2+e_3e_4e_5                           & e_1e_2e_4e_5
        \end{array}\right)
        \end{array}
        $$
If we compute the resolution we get back the same result, thus restating
what is said in \cite[proposition 1.8.]{decker}.

Applying cohomology
algorithms in {\tt Macaulay2} developed by Decker Eisenbud and
Schreyer one may also obtain a full table of the cohomology, which I
write in {\tt Macaulay2} notation:
        \begin{equation} \label{tango-cohomology}
      \vtop{
      \hbox{{\tt {}total:\ 573\ 260\ 92\ 27\ 14\ 7\ 2\ 2\ 7\ 14\ 27\ 92\ 260\ 573}}
      \hbox{{\tt {}\ \ \ -6:\ 573\ 260\ 91\ 14\ \ .\ .\ .\ .\ .\ \ .\ \ .\ \ .\ \ \ .\ \ \ .}}
      \hbox{{\tt {}\ \ \ -5:\ \ \ .\ \ \ .\ \ 1\ 13\ 14\ 7\ 1\ .\ .\ \ .\ \ .\ \ .\ \ \ .\ \ \ .}}
      \hbox{{\tt {}\ \ \ -4:\ \ \ .\ \ \ .\ \ .\ \ .\ \ .\ .\ 1\ .\ .\ \ .\ \ .\ \ .\ \ \ .\ \ \ .}}
      \hbox{{\tt {}\ \ \ -3:\ \ \ .\ \ \ .\ \ .\ \ .\ \ .\ .\ .\ 1\ .\ \ .\ \ .\ \ .\ \ \ .\ \ \ .}}
      \hbox{{\tt {}\ \ \ -2:\ \ \ .\ \ \ .\ \ .\ \ .\ \ .\ .\ .\ 1\ 7\ 14\ 13\ \ 1\ \ \ .\ \ \ .}}
      \hbox{{\tt {}\ \ \ -1:\ \ \ .\ \ \ .\ \ .\ \ .\ \ .\ .\ .\ .\ .\ \ .\ 14\ 91\ 260\ 573}}
      }
        \end{equation}
Reading the table along one antidiagonal gives the list of cohomology
groups of a single twist. Here the list for $T$ starts from the
up-right corner, while starting from a shift to the left means reading
the list for a $(-1)-$twist. \\
One can check that table (\ref{tango-cohomology}) agrees with theorem \ref{bott-cohomology}.

\end{document}